\documentclass[12pt]{article} 
\makeatletter

\newtheorem {theorem}{Theorem}
\newtheorem {corol}{Corollary}
\newtheorem {consec}{Corollary}
\newtheorem {lemma}{Lemma}

\newcommand{\point}{\hspace{-1.75mm}{\bf.\ }}

\newcommand{\btheorem}{\begin{theorem}\point}
\newcommand{\etheorem}{\end{theorem}}

\newcommand{\blemma}{\begin{lemma}\point}
\newcommand{\elemma}{\end{lemma}}

\newcommand{\bpro}{\begin{pro}\point}
\newcommand{\epro}{\end{pro}}

\newcommand{\bcorol}{\begin{corol}\point}
\newcommand{\ecorol}{\end{corol}}

\newcommand{\bnote}{\begin{note}\point}
\newcommand{\enote}{\end{note}}

\newcommand{\bconsec}{\begin{consec}\point}
\newcommand{\econsec}{\end{consec}}

\newcommand{\bdefin}{\begin{defin}\point}
\newcommand{\edefin}{\end{defin}}
\newtheorem {pro}{Proposition}

\begin{document}

\title{A survey of the different types of vector space partitions \thanks{This talk was presented at  Matematiska kollokviet at Department of Mathematics at Link\"oping University on May 19, 2010.}}
\author{Olof Heden}
\maketitle
\begin{abstract}
A {\it vector space partition} is here a collection $\mathcal P$ of subspaces of a finite vector space $V(n,q)$, of dimension $n$ over a finite field with $q$ elements, with the property that every non zero vector is contained in a unique member of $\mathcal P$. Vector space partitions relates to finite projective planes, design theory and error correcting codes.

In the first part of the talk I will discuss some relations between vector space partitions and other branches of mathematics. The other part of the talk contains a survey of known results on the type of a vector space partition, more precisely: the theorem of Beutelspacher and Heden on $\mathrm{T}$-partitions, rather recent results of El-Zanati et al. on the different types that appear in the spaces $V(n,2)$, for $n\leq8$, a result of Heden and Lehmann on vector space partitions and maximal partial spreads including their new necessary condition for the existence of a vector space partition, and furthermore, I will give a theorem of Heden on the length of the tail of a vector space partition.

Finally, I will also give a few historical remarks.
\end{abstract}
\section{Introduction}

We will mainly consider finite dimensional vector spaces over finite fields, and collections of subspaces, covering the whole space, and pairwise intersecting in just the zero vector. Such a configuration $\mathcal P$ of subspaces will be called a {\it vector space partition} of the vector space $V=V(n,q)$, where $n$ is the dimension of $V$ and $q$ the number of elements in the scalar field:
\[
U,U'\in{\mathcal P}\qquad\Longrightarrow\qquad U\cap U'=\{\;\bar 0\;\}
\] 
and
\[
V=\bigcup_{U\in{\mathcal P}}\;U\;.
\]
Let us first give two non trivial, and important examples of vector space partitions.

\bigskip
\noindent
{\bf Example 1.}
Let $q=p^k$ be any power of a prime $p$. Consider the finite field $F=\mathrm{GF}(q^4)$ as a $4$-dimensional vector space $V=V(4,q)$ over the finite field $\mathrm{GF}(q)$. Let $\alpha_1, \alpha_2, \dots,\alpha_k$, where $k=(|F|-1)/(|\mathrm{GF}(q^2)|-1)$, denote a family of coset representatives of the multiplicative group of the subfield $\mathrm{GF}(q^2)$ in  the multiplicative group of $F$. The following family of subspaces of $V$
\[
{\mathcal P}=\{\;\alpha_1\mathrm{GF}(q^2)\;,\;\alpha_2\mathrm{GF}(q^2)\;,\;\dots\;,\;\alpha_k\mathrm{GF}(q^2) \;\}\;
\]
will constitute a vector space partition of $V$.
   
\bigskip
The next construction of a vector space partition is due to Bu \cite{bu} and independently Beutelspacher \cite{beutelspacher}.

\bigskip
\noindent
{\bf Example 2.} Consider the finite field $\mathrm{GF}(q^k)$ as a vector space $W$ over $\mathrm{GF}(q)$, and let $U$ be a subspace of $W$. For each $\alpha\in\mathrm{GF}(q^k)$, we define a subspace $U_{\alpha}$ of $V=W\times U$ by
\[
U_{\alpha}=\{\;(\alpha u,u)\;\mid\;u\in U\;\}\;.
\]
The following set $\mathcal P$ of subspaces to $V$
\[
{\mathcal P}=\{\;U_{\alpha}\;\mid\;\alpha\in\mathrm{GF}(q^k)\;\}\cup\{\;W\times\{\bar 0\}\;\}\;.
\]
will constitute a vector space partition of $V$.
 
\bigskip 
We will say that a vector space partition $\mathcal P$ is of {\it type} 
\[
[d_1^{n_1}d_2^{n_2}\dots d_t^{n_t}]\;,
\]     
if $\mathcal P$ consists of $n_1$ spaces of dimension $d_1$, $n_2$ spaces of dimension $d_2$, etc., where $d_1$, $d_2$, ..., $d_t$ are $t$ distinct non negative integers. So the partition in Example 1 is of type $[2^{q^{2}+1}]$ and the partition in Example 2 is, in case $\dim(W)\neq\dim(U)$, of type $[\dim(W)^1\dim(U)^{|W|}]$.

It might seem to be an easy project to find all possible types of vector space partitions. However, already for finite vector spaces of dimension 8 over a finite field with two elements you get serious problems. In fact this particular case is the ''first'' open case. 

The problem with the possibilities for the different types of vector space partitions were studied in the 70's and 80's by Bu \cite{bu}, Lindstr\"om \cite{bernt}, Beutelspacher \cite{beutelspacher} and Heden \cite{heden}, and during this millennium, and with many contributions by El-Zanati, Seelinger, Sissokho, Spence, and Vanden Eynden, see e.g. \cite{el}, and also by Heden and Lehmann \cite{lehmann}. The subject at issue with this article is to survey these recent results, as well as the results from the 70's and 80's, on this problem.

As a motivation for the study of vector space partitions, we will in the first sections briefly discuss the relation between vector space partitions and projective planes as well as the relation to error correcting codes. 

In the last section we will also give a few historical remarks to vector space partition problems, as well as to group partitions problem.

\section{Vector space partitions and projective planes}\label{sec:2}

A {\it projective plane} consists of lines and points, satisfying the following properties:

\bigskip
1. Any two lines intersect in a unique point,

2. Any two points are contained in a unique line,

3. There are four points such that no line is incident with more than two of them. 

\bigskip
\noindent
Counting arguments show that the number of points will be equal to the number of lines, the integer $q^2+q+1$. The integer $q$ will be called the {\it order} of the plane. Through any point there are exactly $q+1$ lines and every line contains exactly $q+1$ points.

To any projective plane we may associate an {\it affine plane} by deleting one line of the projective plane, the ''line at infinity'' (and the points on that line). What remains will consist of parallel classes of lines. Still, through any two of the remaining points there is a unique line, and further every point is contained in a unique member of every parallel class. Conversely, we may to every given affine plane associate a projective plane by completing with a line at infinity, the points on this line are the distinct parallel classes.   

We will now show a construction due to Andr\'e \cite{andre} of projective planes using vector space partitions. 

Let $V=V(4,q)$ be a $4$-dimensional vector space over a finite field with $q$ elements. Let $\mathcal P$ be any vector space partition of $V$ consisting of solely $2$-dimensional subspaces of $V$:
\[
{\mathcal P}=\{ \;U_1, U_2, \dots, U_t\;\},\quad\hbox{where}\quad \dim(U_i)=2\quad\hbox{for},\quad i=1,2,\dots,t=q^2+1\;.
\] 
(We might as well so far have considered the trivial partition of any $2$-dimensional space into $1$-dimensional spaces.)

By using this partition we first construct an affine plane. The points will be the $q^4$ distinct vectors of $V$. The lines will be the cosets of the spaces in the partition, i.e.,
\[
L_{i,\alpha}=\alpha+U_i\;,\qquad \hbox{for}\qquad \alpha\in\mathrm{GF}(q^4)\;,
\]
(duplications may occur)
and, as distinct cosets of subgroups are disjoint and together cover the whole space, each element $U_i$ of the vector space partition gives a parallel class consisting of $q^4/q^2=q^2$ lines. Further, to show that there is a unique line through any two points $\alpha$ and $\beta$ of the affine plane $\mathrm{GF}(q^4)$, simply, find $U_i$ such that
\[
\alpha-\beta\in U_i\in{\mathcal P}\;.
\] 
Then, the line $L_{i,\beta}=\beta+U_i$ will both contain the point $\beta+0=\beta$ and the point $\beta+(\alpha-\beta)=\alpha$. The verification of the remaining properties of an affine plane are performed in a similar way. By using this affine plane, we may now get a projective plane, as described above, by adjoining a line at infinity. 

The fact is that by using different vector space partitions and the construction of Andr\'e, we get projective planes with different properties. The vector space partition described in Example 1 will give a projective plane that is Desarguessian. But it is rather easy to derive a Non-Desarguessian projective plane using other vector space partitions of the same type.

A vector space partition of the type $[2^{q^2+1}]$ of $V(4,q)$ is called a {\it spread}, or a {\it line spread}, as the members of the spread can be viewed as a family of mutually disjoint lines covering all points in the projective space $\mathrm{PG}(3,q)$ of dimension $3$. 

It must also be remarked that far from every projective plane can be found in this way, see for example the classical book by Dembowski \cite{dembowski}. 

\section{Maximal partial spreads}

One of the most challenging, perhaps most important, but at least a very interesting problem is to make out whether or not there exists a project plane of an order $q$ that is not a power of prime. There are no projective planes of order 6, 10, 14 and for a following known infinite sequence of integer. More precisely: The only general restriction known on the order is given by the Bruck-Ryser-Chowla theorem \cite{bruckryser} which says that if the order $n$ of the projective plane is congruent to 1 or 2 mod 4, it must be the sum of two squares. E.g., 14 is not a sum of two squares. The first undecided case is $q=12$. 

To construct such a plane, we can start with a family of parallel classes of lines, to get a so called {\it partial net}. For such a net it is required that, for every parallel class of lines, every point of the plane is contained in one line of the parallel class. One can get a partial net from an affine plane by deleting all lines in a family of parallel classes of lines. 

It was proved by Bruck \cite{bruck} in 1963, that there is a bound $N(n)$, such that if you have a partial net with more than $N(n)$ parallel classes then you can complete with further parallel classes to get an affine plane of order $n$. More precisely, let 
\[
p(x)=\frac{1}{2}x^4+x^3+x^2+\frac{3}{2}x\;,
\]
and let $d=n-1-t$, where $t$ is the number of parallel classes in a partial net. If
\[
p(d-1)<n\;,
\]
then the partial net can be completed to an affine plane. So, in other words and  with $n=12$: if you find 9 mutually orthogonal latin squares of order 12, then you can always complete with another 2 to get a set of 11 mutually orthogonl latin squares of order 12, sufficiently many to describe a projective plane of order 12. 

These kind of situations can in a natural way be transformed into the field of vector space partition problems, compare the construction of parallel classes of lines from $2$-dimensional subspaces of a $4$-dimensional vector space. 

We define a {\it maximal partial spread} to be a collection $\mathcal S$ of 2-dimensional subspaces of $V=V(4,q)$ with the property that every 2-dimensional subspace of $V$ has a non trivial intersection with at least one member of $\mathcal S$.

The first to study maximal partial spreads was Mesner \cite{mesner}, who in 1967 had his children to choose 2-dimensional spaces, with trivial intersection, ''randomly''. It turned out that if his children found more than a certain bound, they could always complete the lines they had found so far to a full spread. 

Maximal partial spreads have been studied by a numerous of authors as for example Mesner \cite{mesner}, Bruen \cite{bruen}, Bruen and Thas \cite{thas}, Heden \cite{heden93}, Blokhuis \cite{blokhuis}, Heden, Pambianco, Marcugino and Faina \cite{pambianco}, Blokhuis and Metsch \cite{metsch}, Ebert \cite{ebert}, Beutelspacher \cite{albrecht}, G\'acs and Sz\"onyi \cite{gacs}. The best known upper bound for a maximal partial spread is by Blokhuis \cite{blokhuis}:

{\it For every maximal partial spread $\mathcal S$ in $V(4,p)$, where $p$ is a prime number 
\[
|{\mathcal S}|\leq p^2-\frac{p+1}{2}\;.
\]}

Bruen and Thas, Beutelspacher Ebert, and many other researchers, constructed maximal partial spreads of sizes $q^2-q+1$ and $q^2-q+2$; and it was conjectured by Bruen and Thas \cite{thas}, that $q^2-q+2$ is the upper bound for the size of a non trivial maximal partial spread in $V(4,q)$. However, it was proved by Heden \cite{heden00} in 2000 that this conjecture is false, as he found, by using a computer search, a maximal partial spread of size 45 in $V(4,7)$.

We now consider the construction in the Example 2 in a very special, but important case.

\bigskip
\noindent
{\bf Example 3.} We consider the direct product $W\times U$, where $W$ is to be identified with the finite field $\mathrm{GF}(32)$, as a vector space over $\mathrm{GF}(2)$. The space $U$ will be a subspace of dimension 3 of $W$ and by using the construction of Bu, see Example 2, we get a vector space partition of type $[3^{32}5^1]$ of $V(8,2)$. We now continue and partition the subspace $W$ of dimension 5 into one subspace of dimension $3$ and the remaining subspaces of dimension 1. In this way we get a partition of type $[1^{24}3^{33}]$. This will be a partial 3-spread of size 33 in $V(8,2)$. 

If a conjecture of Eisfeld and Storme as well as by Hong and Patel \cite{hong}, was true, then this would be the largest possible size of a partial 3-spread in $V(8,2)$. However, recently, El-Zanati et al \cite{el} found, by using a computer search, a vector space partition of the type $[1^{17}3^{34}]$. This vector space partition will be of some importance below. 

\bigskip
The above example confirms that it is very difficult to find all types of vector space partitions. 

Let us also remark that maximal partial $t$-spreads in $V(n,q)$ have also been extensively studied. Contributions in this study have been given by Beutelspacher \cite{albrecht} in 1980 and during this millenium by Govaerts and Storme \cite{goeverts}.

\section{Vector space partitions and perfect codes}\label{sec:4}

A {\it perfect $e$-error correcting code} is a subset $C$ of a direct product of sets 
\[
C\subseteq {\mathcal A}_1\times{\mathcal A}_2\times\dots\times{\mathcal A}_t
\]
such that any possible word $x$ of this direct product differs in at most $e$ coordinate positions from a unique word of $C$. Error correcting codes in general, not just perfect error correcting codes, is a well studied subject, originating in the late 40's during the development of computers, but is also of great importance, in fact most important, in connection with information transmission.

Herzog and Sch\"onheim \cite{herzog} observed in 1972 that vector space partitions can be used to construct perfect 1-error correcting codes. Take any vector space partition
\[
{\mathcal P}=\{\;U_1, U_2, \dots, U_t\;\}
\]
of the space $V=V(n,q)$. Consider the map $\varphi$ from the direct product of the spaces of $\mathcal P$,
\[
\varphi:\;U_1\times U_2\times\dots \times U_t\qquad \longrightarrow\qquad V
\]
defined by
\[
(u_1,u_2,\dots,u_t)\quad\mapsto\quad u_1+u_2+\dots+u_t\;.
\]
The kernel of this map, i.e.,
\[
\ker(\varphi)=\{\;(u_1,u_2,\dots,u_t)\;\mid\;u_1+u_2+\dots+u_t=0\;\}
\]
will be a perfect 1-error correcting code. In case the spaces in the vector space partition $\mathcal P$ are not all of the same size, or equivalently, not of the same dimension, then the codes are called {\it mixed perfect codes}. 

Those acquainted with the ''ordinary'' Hamming code may recognize the above construction as a generalization of the well known construction of a Hamming code as the null space of a matrix $\bf H$. For example, consider
\[
{\bf H}=\left(\begin{array}{ccccccc}
0&0&0&1&1&1&1\\
0&1&1&0&0&1&1\\
1&0&1&0&1&0&1
\end{array}\right)
\]
In fact the Hamming code we get as null space of this matrix, can be obtained by the vector space partition consisting of the following subspaces of $V(3,2)$:
\[
{\mathcal P}=\{\;U_1=\{(0,0,0),(0,0,1)\}\,,\;U_2=\{(0,0,0),(0,1,0)\}\,,\;\dots,\;U_7=\{(0,0,0),(1,1,1)\}\;\}
\]

All perfect codes, constructed in this way, will actually be linear codes, that is, vector spaces over the field $\mathrm{GF}(q)$.

Finally, as was observed by Herzog and Sch\"onheim \cite{herzog}, every linear perfect 1-error correcting code originates from a vector space partition, as described above.

\section{On the types of vector space partitions}

As already mentioned, and perhaps also realized from the Example 3 above, to find all possible types of vector space partitions is a difficult task. It contains both finding new constructions and proving necessary conditions for a certain type to exist. So far no general necessary and sufficient conditions have been found. 

\subsection{Necessary conditions}

We will in this subsection always assume that if $[d_1^{n_1}d_2^{n_2}\dots d_k^{n_k}]$ is a type of a vector space partition then $d_1<d_2<\dots<d_k$.

As any vector is contained in a unique space of the partition, the following so called {\it packing condition} for a vector space partition of type $[d_1^{n_1}d_2^{n_2}\dots d_k^{n_k}]$ in $V=V(n,q)$ to exist must be true:
\[
n_1(q^{d_1}-1)+n_2(q^{d_2}-1)+\dots+n_k(q^{d_k}-1)=q^n-1\;.
\]
As for any two members $U$ and $W$ of a vector space partition $\mathcal P$ of $V=V(n,q)$ we have that $\dim(\mathrm{span}\{U\cup W\})=\dim(U)+\dim(W)$ we get, as first observed by Bu \cite{bu}, for any $i$ and $j$:
\[
d_i+d_j\leq n\;.
\] 
This condition will below be called the {\it dimension condition}. Example 2 shows that we can have equality above. 

Spera \cite{spera} observed that
from the packing condition follows that,
\[
\frac{q^n-1}{q^{d_k}-1}\leq \sum_{i=1}^kn_i=|{\mathcal P}|\leq \frac{q^n-1}{q^{d_1}-1}\;.
\]
Heden and Lehmann \cite{lehmann} improved this bound for the number of spaces in a vector space partition: If we are not in the case of Example 2 of Section 1, then
\[
|{\mathcal P}|\geq q^{d_k}+q^{d_{k-1}}+1\;.
\]

We used the so called second packing condition developed from the fact that if we intersect all spaces of a vector space partition $\mathcal P$ with a hyperplane $H$, then we get a vector space partition ${\mathcal P}_H$. To describe the second packing condition it will be convenient to use the following notation:

{\it A $(m_k,m_{k-1},\dots,m_2,m_1)$-partition is the same as a $[1^{m_1}2^{m_2}\dots k^{m_k}]$-partition, (where we allow some of the non negative integer exponents to be zero).} 

A hyperplane $H$ will be of type $b=(b_k,\dots,b_2,b_1)$ if $H$ contains $b_i$ of the subspaces of dimension $i$ of $\mathcal P$. Let $s_b$ denote the number of hyperplanes of type $b$. The {\it second packing condition}, derived by Heden and Lehmann \cite{lehmann}, is
\[
s_{b}\neq0\qquad\Longrightarrow\qquad\sum_{d=1}^{k} b_{d}q^{d}=\sum_{d=1}^{k}m_{d} - 1\;.
\] 
Let $B$ denote the family of all feasible solutions to the diophantine equation above. By double counting incidences, Heden and Lehmann \cite{lehmann} proved the following necessary conditions:

{\it For any $1\leq d,d'\leq n-2$,
\begin{eqnarray}
\sum_{b\in B}b_{d}s_{b} & = & m_{d}\\
\sum_{b\in B}{{b_{d}}\choose{2}}s_{b} & = & {{m_{d}}\choose{2}}\\
\sum_{b\in B}b_{d}b_{d'}s_{b} & = & m_{d}m_{d'}\;.\label{eq:sum-sb-i-j}
\end{eqnarray}}

Using these necessary conditions, that we shall call the {\it hyperplane conditions}, Heden and Lehmann \cite{lehmann} derived the following 

{\it Let $V=V(2t,q)$ and assume that
$V$ has a partition of type $(m_t,\ldots,m_1)$, with $m_t=q^t+1-a$. Let $d<t$, such that $m_d>0$. 

If \begin{equation}
m_d<\frac{q^{t}-1}{q^{t-d}-1}\;,\end{equation}
 then \begin{equation}
a\geq m_d-\mathrm{R}_q(t,d,m_d)\;.\end{equation}
where
\[
\mathrm{R}_q(t,d,m):=m(m-1)\frac{\frac{1}{2}(q^{2t-2d}-1)+1-q^{t-d}}{q^{t}-1-m(q^{t-d}-1)}\;.
\]}

We now evaluate this bound in a particular case. We consider as above $V(2t,q)$. We can easily, like in Example 1 of Section 1, find a spread consisting of $q^t+1$ spaces of dimension $t$. We can then substitute a number $a$ of these spaces, by vector spaces partitions of them, consisting of spaces of lower dimension. If $t/2<d<t$ then, in this way, we cannot get a vector space partition with more than $a$ subspaces of dimension $d$. 

The problem Heden and Lehmann \cite{lehmann} considered was whether it is possible to get more than $a$ spaces of dimension $d$, still with $q^t+1-a$ of spaces of dimension $t$, (but then, as spaces of dimension $t$, taking other spaces than those from a complete spread.) However, they proved that if the number of spaces of dimension $d$ is small then this is not possible, more precisely: 

{\it Let $d=t-k$ and let $m_d$ be the number of spaces of dimension $d$. If
\[
m_{d}\leq \sqrt{2}q^{(t-2k)/2}\;,
\]
then $a\geq m_{d}$.
}  

We close this subsection by discussing the length of the tail. The {\it tail} of a vector space partition is the set of spaces of lowest dimension. The {\it length of the tail} is the size of the tail. The size of the tail in the vector space partition in Example 3, which actually will be of importance below, is $17$.

By using the connection with perfect codes, as described in Section \ref{sec:4}, Heden \cite{h} proved the following bounds for the length of the tail of a vector space partition.

{\it For every vector space partition of type $[d_1^{n_1}d_2^{n_2}\dots d_k^{n_k}]$

\bigskip
\begin{tabular}{ll}
(i)&if $q^{d_2-d_1}$ does not divide $n_1$ and if $d_2< 2d_1$, then $n_1\geq q^{d_1}+1$.\\
(ii)& if $q^{d_2-d_1}$ does not divide $n_1$ and $d_2\geq 2d_1$, then either\cr & $d_1$ divides $d_2$ and $n_1=(q^{d_2}-1)/(q^{d_1}-1)$ or $n_1> 2 q^{d_2-d_1}$\\
(iii)& if $q^{d_2-d_1}$ divides $n_1$ and $d_2<2d_1$ then $n_1\geq q^{d_2}-q^{d_1}+q^{d_2-d_1}$.\\ 
(iv)& if $q^{d_2-d_1}$ divides $n_1$ and $d_2\geq2d_1$ then $n_1\geq q^{d_2}$.
\end{tabular}}

\bigskip
It was shown in \cite{h} that the bounds in {\it (i), (iii)} and {\it (iv)} are tight.

\subsection{$T$-partitions}

A vector space partition $\mathcal P$ is a {\it $T$-partition} if
\[
T=\{\;\dim(W)\;\mid\;W\in{\mathcal P}\;\}\;.
\]

The problem is to find conditions on a given set $T$ of positive integers, that guarantees the existence of a $T$-partition. We will always assume that
\[
T=\{\;t_1, t_2, \dots, t_k\;\} \qquad\hbox{where}\qquad t_1<t_2<\dots<t_k\;.
\]

In the case $t_1=1$, it is trivial to find a $T$-partition of any space $V=V(n,q)$ where $t_k+t_{k-1}\leq n$. Simply, consider $V$ as a direct product of $W$ and $U$ where $d_2=\dim(W)=t_k$ and $d_1=\dim(U)=n-\dim(W)$, and so using the construction of Example 2 of a vector space partition, we get a partition of type $[d_1^{q^{d_2}}d_2^1]$. We now take $k-1$ of the subspaces of dimension $d_1$ and partition each of them into one space of dimension $t_i$ and the remaining subspaces of dimension $1$, respectively. We note that $k\leq d_2$ implies that $k<q^{d_2}$, so this type of $T$-partition will exist. 

So the non trivial case is when $t_1\geq2$. 

It was Beutelspacher \cite{beutelspacher} who in 1978 introduced the concept $T$-partition. He found a result on the existence of $T$-partitions that relates to the Frobenius number.

Let $A=\{\;a_1,a_2,\dots,a_k\;\}$ be a set of positive integers and assume that the greatest common divisors of these numbers is 1. The greatest integer $n$ that cannot be written as a linear combination of these integers
\[
n=x_1a_1+n_2x_2+\dots+n_ka_k\;
\]
for some non negative integers $x_1$, $x_2$, ..., $x_k$, is the {\it Frobenius number} $g(A)$. It was proved by Selmer \cite{selmer} that
\[
g(A)\leq 2a_1\lfloor \frac{a_k}{k}\rfloor-a_1\;,
\]
where $a_1$ is the smallest, and $a_k$ the largest among the integers in $A$. 

Using the above result of Selmer, Beutelspacher \cite{beutelspacher} proved the following

{\it
Consider the vector space $V=V(n,q)$. For
\[
n>2t_1\lfloor \frac{t_k}{d\cdot t_k}\rfloor+t_2+\dots+t_k\;, 
\]
$V$ has a $T$-partition if and only if $gcd(T)$ divides $n$.
}

It must be remarked that from the packing condition follows that if $V(n,q)$ admits a $T$-partition, then $\gcd(T)$ divides $n$.

Further, Beutelspacher \cite{beutelspacher} proved the next theorem for $T$-partitions, in the case $t_1=2$, and later, Heden \cite{heden} proved the result in its full generality:

{\it
The space $V(2t,q)$ admits a $T=\{t_1<t_2<\dots<t_k=t\}$-partition.
}

\subsection{An enumeration of the different types of vector space partitions in $V(n,2)$, for $n\leq 7$.}

This enumeration was completed by El-Zanati et al in \cite{esssv}, although the investigation of almost all cases were ruled out already by Heden \cite{He3}. It was shown, partially by the use of a computer search for one particular vector space partition, that the following necessary conditions: the packing condition, the dimension condition and the tail condition also are sufficient in the case $V(n\leq7,2)$. For  e.g. $n\leq5$ we thus get the following enumeration:
\[
\begin{array}{c|l}
n&\hbox{different types of vector space partitions in $V(n,2)$}\\
\hline
1&[1^1],\\
2&[1^3],\\
&[2^1],\\
3&[1^7],\;[1^42^1],\\
&[3^1],\\
4&[1^{15}],\;[1^{12}2^{1}],\;[1^92^2],\;[1^{6}2^{3}],\;[1^32^4],\\
&[2^5],\\
&[1^{8}3^1],\\
&[4^1],\\
5&[1^{31}],\;[1^{28}2^1],\;[1^{25}2^2],\;[1^{22}2^3],\;[1^{19}2^4]\;[1^{16}2^5],\;[1^{13}2^6],\;[1^{10}2^7],\;[1^72^8],\;[1^42^9],\\
&[1^{24}3^1],\;[1^{21}2^13^1],\;[1^{18}2^23^1],\;[1^{15}2^33^1],\;[1^{12}2^43^1],\;[1^92^53^1],\;[1^62^63^1],\;[1^32^73^1],\;[2^83^1],\\
&[1^{16}4^1],\\
&[5^1].
\end{array}
\] 

Together with Heden, El-Zanati et al \cite{hedenpapa} investigated the case $n=8$ and $q=2$ and vector space partitions consisting of spaces of dimension at least equal to $2$. It turned out that the packing condition, dimension condition and the tail condition, with just one exception, were both necessary and sufficient in this case. The exceptional case, that could not be excluded by these three conditions, was the existence of a vector space partition of type $[2^63^64^{13}]$, which however could be excluded by the use of the hyperplane conditions, specialized to this case. The non existence of that vector space partition thus also follows from the theorem of Heden and Lehmann \cite{lehmann}.

\subsection{Are there any conditions that are both necessary and sufficient?}

It was shown by Heden and Lehmann \cite{lehmann}, with the use of examples, that the four necessary conditions, the packing-, the dimension-, the tail-, and the hyperplane condition are not sufficient for the existence of a vector space partition. We now give the very few and very special instances, but general according to dimension and size of the scalar field, in which we have a complete picture of the different types.

As can easily be proved by elementary arguments, $q^d-1$ divides $q^n-1$, for a prime power $q$, if and only if $d$ divides $n$. Hence, by first packing condition, a necessary condition for a vector space partition of type $[d^m]$ to exist in $V(n,q)$ is that $d$ divides $n$. This is also sufficient, as then $\mathrm{GF}(q^n)$ has a subfield $\mathrm{GF}(q^d)$ and we can easily find a vector space partition, using the same construction as in Example 1, of type $[d^m]$ where $m=(q^n-1)/(q^d-1)$.

However, as soon as we leave the above simple case, the situation gets complicated, even in the case of just two different dimensions appearing in the vector space partition. As was indicated by Example 3, the case $[1^{n_1}d^{n_2}]$ is still far from completely investigated.

In his thesis, Heden \cite{heden09} proved the following

{\it
The packing condition, dimension condition together with the condition 
\[
\dim(U_i)\geq c\;,\qquad\hbox{for $i\leq q$}\;,
\] 
are necessary and sufficient for the existence of a vector space partition $U_1$, $U_2$, ..., $U_k$ of $V(n,q)$, if $\dim(U_{q+1})=\dim(U_{q+2})=\ldots=\dim(U_k)=c$.
}

This theorem is a generalization of a theorem of Lindstr\"om \cite{bernt}, who proved the theorem in the case all but one of the subspaces have dimension $c$.

\section{Some historical remarks}

The first who investigated these kind of problems was George Abram Miller. In 1906 he published a paper \cite{miller} in which he proved that if an abelian group $G$ admits a partition into subgroups, then all elements of $G$ must have order $p$, for some prime number $p$. 

The proof idea is simple. Assume that there are elements $h_i$ and $h_j$ of order $p$ and $q$ respectively, where $p$ is a prime, and such that $h_i$ and $h_j$ are elements of different subgroups $H_i$ and $H_j$ of $G$ in the partition of $G$. Then, $h_i+h_j=h_k$ is an element in a third group $H_k$ of the partition. The sum of $h_k$, a $p$ numbers of times, will give
\[
p\cdot h_k=(h_i+h_j)+(h_i+h_j)+\dots+(h_i+h_j)=p\cdot h_j\;,
\]
which is an element of both $H_k$ and $H_j$, and hence zero. Miller also found a group partition of a group with $p^2$ elements into subgroups, each with $p$ elements.

Twenty years later John Wesley Young \cite{young}, who had made his master thesis under the supervision of Miller, studied the problem of the partition of infinite groups into subgroups. A Russian, Kontorovich \cite{kontorovich} studied and published results in 1939 and 1940 about partitions of a group $G$ having a special property: $G=HK=KH$ for any two members $H$ and $K$ in the partition. 

Let us give a simple example of a partition of a non abelian group. The group ${\mathcal S}_3$ has a partition
\[
{\mathcal S}_3=\{\mathrm{id.},\;(1\;2)\}\cup\{\mathrm{id.},\;(1\;3)\}\cup\{\mathrm{id.},\;(2\;3)\}\cup\{\mathrm{id.},\;(1\;2\;3),\;(1\;3\;2)\}\;.
\]

Later, the main problem studied, mainly by Reinhold Baer \cite{baer} and his student Otto Kegel \cite{kegel} and by Michio Suzuki \cite{szuzuki}, was to classify those non abelian groups that admits a partition into subgroups. 

A nice survey of this search for a classification was given by Zappa \cite{zappa} in 2003. Quoting Zappa, a group $G$ has a non trivial partition if and only if it satisfies one of the following conditions:

\begin{enumerate}
\item $G$ is a $p$-group with $H_P(G)\neq G$ and $|G|>p$;
\item $G$ is a Frobenius group;
\item $G$ is a group of Hughes-Thompson type;
\item $G$ is isomorphic with $\mathrm{PGL}(2,p^h)$, $p$ being an odd prime;
\item $G$ is isomorphic with $\mathrm{PSL}(2,p^h)$, $p$ being a prime;
\item $G$ is isomorphic with a Suzuki group $G(q)$, $q=2^h$, $h>1$.
\end{enumerate} 

We must remark that $H_p(G)$ is the so called {\it Hughes subgroup}, i.e., the subgroup of $G$ generated by those elements of $G$ that have not the order $p$. A {\it Frobenius group} is a transitive permutation group on a finite set, such that no non-trivial element fixes more than one point and some non-trivial element fixes a point. The group ${\mathcal S}_3$ is an example of a Frobenius group.


\begin{thebibliography}{12}

\bibitem{andre}
J. Andr\'e, \"Uber nicht-Desarguesschen Ebenen mit transitiven Translationsgruppe, Math. Zeitschr, 60(1954)156--186.

\bibitem{baer}
R. Baer, Partitionen endlicher Gruppen, Math Zeitschr., 75(1960/61)333--372.

\bibitem{beutelspacher}
A. Beutelspacher, Partitions of finite vector spaces: an application of the Frobenius number in geometry, Arch. Math., 31(1978)202--208.

\bibitem{albrecht}
A. Beutelspacher, Blocking sets and partial spreads in finite projective spaces, Geometriae Dedicata, 9(1980)425--449.

\bibitem{blokhuis}
A.Blokhuis, On the size of a blocking set in $PG(2,p)$, Combinatorica, 4(1)(1994)111--114.

\bibitem{metsch}
A. Blokhuis, K. Metsch, On the size of a maximal partial spread, Designs, Codes and Cryptography, 3(1993)187--191

\bibitem{bruck}
R. H. Bruck, Finite nets. II. Uniqueness and imbedding. Pacific J. Math., 13(1963)421-457.

\bibitem{bruckryser}
R. H. Bruck, H. J. Ryser, The non existence of certain finite projective planes, Canad. J. Math., 1(1949)88--93.

\bibitem{bruen}
A. A. Bruen, Partial spreads and replaceable nets, Canad. J. Math., 23(1971)381--392.

\bibitem{thas}
A. A. Bruen, J Thas, Partial spreads, packings and Hermitian manifolds in $PG(3,q)$, Math. Zeitschr., 151 (1976), 207--214. 

\bibitem{bu} T.\ Bu, Partitions of a vector space, {Discrete
Math.}, {31}(1980)79--83.

\bibitem{dembowski}
P. Dembowski, Finite Geometries, Springer, 1997.

\bibitem{ebert}
G. Ebert, Maximal strictly partial spreads, Canad. J. Math, 30(1978)483--489.

\bibitem{esssv}
S.\ El-Zanati, G.\ Seelinger, P.\ Sissokho, L.\ Spence, and C.\ Vanden Eynden, 
On partitions of finite vector spaces of small dimension over $GF(2)$, {Discrete Mathematics}, {309}(2009)4727--4735.

\bibitem{el}
S.\ El-Zanati, G.\ Seelinger, P.\ Sissokho, L.\ Spence, and C.\ Vanden Eynden, 
The maximum size of a partial $3$-spread in a finite vector space over $\mathrm{GF}(2)$, Designs, Codes and Cryptography, 54(2010)101--107.

\bibitem{hedenpapa}
S. El-Zanati, O. Heden, G. Seelinger, P. Sissokho, L. Spence,
C.~Vanden~Eynden, Partitions of the $8$-dimensional vector space over $GF(2)$, Journal of Combinatorial Designs, 18(2010)462--474.

\bibitem{gacs}
A. G\'acs, T. Sz\"onyi, On maximal partial spreads in $\mathrm{PG}(n,q)$, Des. Codes and Cryptogr, 29(2003)123--129. 

\bibitem{goeverts} P. Govaerts, L. Storme, On a particular class
of minihypers and its applications. I. The result for general $q$,
Designs, Codes and Cryptography, 28(2003)51--63.

\bibitem{heden}
O. Heden, The Frobenius number and partitions of a finite vector space, Arch. Math., 42(1984)185--192.

\bibitem{He3} O.\ Heden, Partitions of finite abelian groups,
{ Europ. J. Combin.}, {7}(1986)11--25.

\bibitem{heden93}
O. Heden, Maximal partial spreads and the modular $n$-queen problem, Discrete mathematics, 120(1993), 75--91.

\bibitem{heden00}
O. Heden, A maximal partial spread of size 45 in $\mathrm{PG}(3,7)$, Designs, Codes and Cryptography, 22(2000)331--334.

\bibitem{heden09}
O. Heden, Necessary and sufficient conditions for the existence of a class of partitions of a finite vector space, Designs, Codes and Cryptography, 53(2009)69--73.

\bibitem{h}
O. Heden, On the length of the tail of a vector space partition, Discrete mathematics, 309(2009)6196--6180.

\bibitem{pambianco}
O. Heden, G. Faina, S. Marcugini, F. Pambianco, The maximal size of a maximal partial spread in $PG(3,9)$, submitted (2006).

\bibitem{lehmann}
O. Heden, J. Lehmann, Some necessary conditions for vector space partitions, submitted.

\bibitem{herzog}
M. Herzog and J. Sch\"onheim, Group partition, factorization and the vector covering problem, Canad. Math. Bull., 15(2)(1972) 207--214. 

\bibitem{hong}
S. Hong, A. Patel, A general class of maximal codes for computer applications, IEEE Trans. Comput., C-21,(1972)1322--1331.

\bibitem{kegel}
O. Kegel, Nich-einfache Partitionen endlicher Gruppen, Arch. Math., 12(1961)170--175. 

\bibitem{kontorovich}
P. Kontorovich, On the representation of a group as a direct product of its subgroups II, Mat. Sb, 7(1940)27--33 (in Russia).

\bibitem{bernt}
B. Lindstr\"om, Group partitions and mixed perfect codes, Canad. Math. Bull., 18(1)(1975)57--60.

\bibitem{mesner}
D. Mesner, Sets of disjoint lines in $PG(3,q)$, Canad. J. Math., 19(1967)273--280.

\bibitem{miller}
G. A. Miller, Groups in which all operators are contained in  series of subgroups such that any two have only identity in common, Bull. Amer. Math., 12(1905-1906)446--449.

\bibitem{selmer}
E. S. Selmer, On the linear diophantine problem of Frobenius, J. reine angew. Math., 293(1977)1--17.

\bibitem{spera} 
A. Spera, On partitions of finite vector spaces, arXiv:0902.3075v1, 18 Feb 2009.

\bibitem{szuzuki}
M. Suzuki, On the finite groups with a complete partition, J. Math. Soc. Japan, 2(1950)165--185. 
\bibitem{young}
J. W. Young, On the partition of a group and the resulting classification, Bull. Amer. Math. Soc., 35(1927)453--461.

\bibitem{zappa}
G. Zappa, Partitions and other coverings of finite groups, Illinois Journal of Mathematics, 47(2003)571--580.
\end{thebibliography}
\end{document}